\documentclass[12pt]{article}

\usepackage{amssymb,amsmath,latexsym,times,verbatim}
\usepackage{bm}
\usepackage{amsfonts}
\usepackage{mathrsfs}
\usepackage{graphicx,epsfig,epstopdf,color}

\newcommand{\ba}{\begin{array}}
\newcommand{\ea}{\end{array}}

\newcommand{\R}{{\mathbb R}}

\newcommand{\C}{{\mathbb C}}

\newcommand{\PROOF}{{\bf Proof:} }

\renewcommand{\epsilon}{\varepsilon}

\newcommand{\qed}{\hfill \rule{3mm}{3mm}}

\newtheorem{theorem}{Theorem}[section]

  {\begin{trivlist}\item[]\textbf{Proof#1 }}%
  {\hspace*{\fill}$\rule{.3\baselineskip}{.35\baselineskip}$\end{trivlist}}

\def\eps{\epsilon}

\begin{document}

\title{On the mathematical description of time--dependent surface water waves}
\author{\renewcommand{\thefootnote}{\arabic{footnote}} Wolf-Patrick D\"ull$^{1}$}

\footnotetext[1]{Institut f\"ur Analysis, Dynamik und  Modellierung, 
       Universit\"at Stuttgart, Pfaffenwaldring 57, 70569 Stuttgart, Germany,
duell@mathematik.uni-stuttgart.de}
\maketitle

\begin{abstract}
This article provides a survey on some main results and recent developments in the mathematical theory of water waves. More precisely, we briefly discuss the mathematical modeling of water waves and then we give an overview of local
and global well--posedness results for the model equations. Moreover, we present reduced models in various parameter regimes for the approximate description of the motion of typical wave profiles and discuss the mathematically rigorous justification of the validity of these models.   
\\[3mm]
{\bf Keywords\;} water wave problem, well--posedness, reduced models, approximation
\\[3mm]
{\bf Mathematics Subject Classification (2010)\;} 76B15, 35Q35, 35Q53, 35Q55

\end{abstract}
%\tableofcontents{}

\section{Introduction}
%A better prediction of the dynamics of surface water waves would be of %great practical and theoretical relevance. One need only think of the %development of efficient early warning systems for tsunamis. 
The mathematical analysis of water waves has a long history of more than 300 years. Over the last two decades the further development of the mathematical theory of water waves has attracted new and increasing interest leading to significant progress.
Nevertheless, there are still many challenging open problems being both relevant for practical applications and stimulating for the development of new subtle mathematical concepts. This article provides
 a survey on some main results and recent developments in the modern mathematical theory of water waves with special emphasis on the local and global well--posedness theory of the model equations and on the justification of reduced models. Since the existing literature is very extensive, a short survey article on the water--wave theory cannot be exhaustive. Hence, we refer  for further information to the historical articles \cite{Da03, Ck04, Ck05}, 
the surveys \cite{SW02b, CW07}, the books \cite{St92, LM76, De94, J97, C11, L13} and the references therein.

Mathematically, water waves are modeled by the Navier--Stokes or Euler's
equations with some free unknown top surface. Usually viscosity and compressibility can be neglected such that the incompressible Euler's equations are chosen. The resulting mathematical model is called the water wave problem. 

In many practically relevant problems, it is sufficient to study the water wave problem only in two dimensions -- one for the height and one for the direction of the propagation of the wave -- to get a realistic description of the problems. For example, if an earthquake happens in a crack along a continent it creates almost parallel initial conditions for tsunamis. 

We now discuss the water wave problem in more detail. For simplicity we restrict ourselves to the 2--dimensional case with finite depth of water and a flat bottom. The modifications for uneven bottoms, for infinite depth of water and for the 3--dimensional case are straightforward.

\subsection{The water wave problem in Eulerian coordinates}

\begin{figure}[htbp]
\centering
\epsfig{file=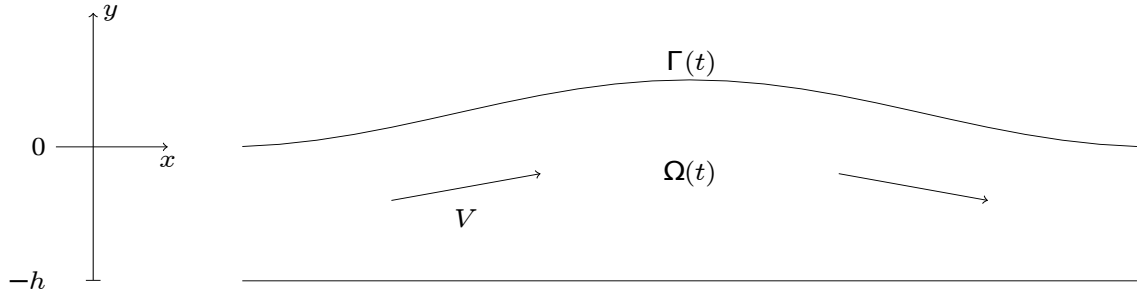,width=15cm,height=5cm,angle=0}
\vspace*{-1.1cm}
\caption{{\small The 2--D water wave problem with finite depth.} \label{fig1} }
\end{figure}

The 2--D water wave problem for waves over an incompressible, inviscid fluid in an infinitely long canal of finite depth under the influence of gravity and surface tension has in Eulerian coordinates the following form: The fluid fills a domain $\Omega(t)=\{(x,y) \in \R^2 \,|\, x \in \R, -h < y < \eta(x,t)\}$ in between the bottom  $B = \{(x,y) \in \R^2 \,|\, x \in \R,\, y = -h\}$ and the free top surface $\Gamma(t) = \{(x,y) \in \R^2 \,|\, x \in \R,\, y = \eta(x,t)\}$; see Figure \ref{fig1}.

\vspace*{0.1cm}

The velocity field $V=(v_1, v_2)$ of the fluid is governed by the incompressible Euler's equations 
\begin{align}
V_t + (V\cdot\nabla)V & =  -\nabla p + g \left(\!\begin{array}{c} 0 \\ -1 \end{array} \!\right) \qquad \text{in}\,\, \Omega(t), \label{euler} \\[2mm] 
\nabla \cdot V & =  0 \qquad \text{in}\,\, \Omega(t), \label{incompr}
\end{align}
where $p$ is the pressure and $g$ the constant of gravity. 

Assuming that fluid particles on the top surface remain on the top
surface, that the pressure at the top surface is determined by the
Laplace--Young jump condition and that the bottom is impermeable yields the boundary conditions
\begin{align}
\eta_t & =  V \cdot \left(\!\begin{array}{c} -\eta_x \\ 1 \end{array} \!\right) \qquad \text{at}\,\, \Gamma(t), \label{surf} \\[2mm]  
p & =  -bgh^2 \kappa \qquad \text{at}\,\, \Gamma(t), \label{press} \\[2mm]  
v_2 & =  0   \qquad \text{at}\,\, B \label{bot},
\end{align}
where $b \geq 0$ is the Bond number, which is proportional to the strength of the surface tension, and $\kappa$ is the curvature of $\Gamma(t)$.

If the flow is additionally assumed to be irrotational, the above
system can be reduced to a system defined on $\Gamma(t)$. Due to the
irrotationality of the motion there exists a velocity potential $\phi$
with $V=\nabla \phi$, which is harmonic in $\Omega(t)$ with vanishing
normal derivative at $B$. Moreover, the motion of the vertical component of the velocity is uniquely determined
by the horizontal one, i.e., there exists an operator 
$ \mathcal{K} = \mathcal{K}(\eta) $ such that
\begin{align} \label{Kop}
\phi_y &= \mathcal{K} \phi_x.
\end{align}
The operator $ \mathcal{K}$ is related to the Dirichlet to Neumann operator.

By using the potential $\phi$, the system (\ref{euler})--(\ref{bot})
can be reduced to  
\begin{align} 
\eta_t & =  V \cdot \left(\!\begin{array}{c} -\eta_x \\ 1 \end{array} \!\right) \qquad \mathrm{at}\,\, \Gamma(t), \label{isurf2} \\[2mm]  
\phi_t &= -\textstyle{\frac{1}{2}} ((\phi_x)^{2} + 
(\mathcal{K} \phi_x)^{2}) - g \eta + b g h^2 \left( \frac{\eta_x}{\sqrt{1+\eta_x^2}}\right)_{x} \qquad \text{at}\,\, \Gamma(t) \label{ipot}
\end{align}
or to
\begin{align} 
\eta_t & =  \mathcal{K} v_1 -v_1\eta_x 
 \qquad \text{at}\,\, \Gamma(t), \label{surf2} \\[2mm]  
(v_1)_t &= - g \eta_x -\textstyle{\frac{1}{2}} (v_1^{2} + 
(\mathcal{K} v_1)^{2})_x + b g h^2 \left( \frac{\eta_x}{\sqrt{1+\eta_x^2}}\right)_{xx} \qquad \text{at}\,\, \Gamma(t). \label{pot}
\end{align}
Alternatively, by taking the trace of equation \eqref{euler} at $\Gamma(t)$ 
and using \eqref{surf}, \eqref{press} and \eqref{Kop}, the system \eqref{euler}--\eqref{bot} can be reduced to 
\begin{align} 
\eta_t & =  \mathcal{K} v_1 -v_1\eta_x 
 \qquad \text{at}\,\, \Gamma(t), \label{surf3} \\[2mm]  
(v_1)_t &= - a \eta_x  -v_1 (v_1)_x + b g h^2 \left( \frac{\eta_x}{\sqrt{1+\eta_x^2}}\right)_{xx} \qquad \text{at}\,\, \Gamma(t), \label{pot3}
\end{align}
with
\begin{align} 
a= - p_y = g + (\mathcal{K} v_1)_t + v_1 (\mathcal{K} v_1)_x\,. \label{acoeff} 
\end{align}
The system \eqref{surf3}--\eqref{acoeff} even holds for flows with vorticity if $
\mathcal{K} v_1$ is replaced by $v_2$. However, in this case, the evolution equation for the vorticity or a related quantity on the entire domain $\Omega(t)$ is needed to close the system.

\subsection{Alternative coordinate systems}

Choosing Eulerian coordinates to formulate the equations for the motion of water waves is natural for describing many physical experiments. However, for surface waves which cannot be described by a graph using Eulerian coordinates is not so convenient. But there are also alternative coordinate systems which yield appropriate frameworks for formulating the water wave problem. Each of these coordinate systems has its own advantages concerning applicability and mathematical structure of the resulting equations of the water wave problem. Hence, depending on the problem one intends to solve, one has to find out which coordinate system is the most adapted one. In the following, we give an overview of the most important alternative coordinate systems which are used to formulate the water wave problem.
\bigskip

The second classical formulation of the water wave problem bases on Lagrangian coordinates. In these coordinates, the free top surface is parametrized as
$\Gamma(t) = \{X(\alpha,t)\,|\, \alpha \in \R\}$,
where $X(\alpha,t)$ are the Eulerian coordinates of the fluid particle labeled by $\alpha$. This parametrization is well--defined due to the above assumption that fluid particles on the top surface remain on the top surface.
We have
\begin{equation}
{X}_t(\alpha,t) = V({X}(\alpha,t),t)
\end{equation}
such that \eqref{euler} implies
\begin{equation} \label{lagr1}
{X}_{tt}(\alpha,t)=-\nabla p ({X}(\alpha,t),t) + g \left(\!\begin{array}{c} 0 \\ -1 \end{array} \!\right)\,.
\end{equation}
One advantage of this equation is that the transport term $(V\cdot\nabla)V$ has disappeared.

Assuming \eqref{incompr}, \eqref{press}, \eqref{bot} and the irrotationality of the flow, the dynamics of the fluid is again completely determined by the evolution of the free surface. In the literature, there are various evolutionary systems for the free surface derived from  \eqref{lagr1}; see, for example, \cite{Na74, Y82, Y83, Cr85, Wu97, Wu99, SW00, SW02a}.
\bigskip

Another formulation of the 2--D water wave problem is the so-called arc length formulation, where the free
top surface is parametrized by arc length. This formulation was introduced in the theory of the 2--D water wave problem by Ambrose and Masmoudi \cite{A03, AM05}, influenced by the numerical work of Hou, Lowengrub and Shelley \cite{HLS94, HLS95}. Let $P(t): \R \rightarrow \Gamma(t),\,\alpha \mapsto P(\alpha,t) = (x(\alpha,t), y(\alpha,t))$ be a parametrization of $\Gamma(t)$ by arc length. Hence, we have $(x_{\alpha}^2+y_{\alpha}^2)^{1/2}=1$. 
Let $U$ and $T$ be the normal and the tangential velocity field of the free top surface measured in the coordinates of the arc length parametrization, that means that
\begin{equation} \label{xy}
(x,y)_{t}(\alpha,t) = U(\alpha,t) \hat{n}(\alpha,t) + T(\alpha,t) \hat{t}(\alpha,t),
\end{equation}
where $\hat{n} =(-\sin(\theta), \cos(\theta))$ and $\hat{t}=(\cos(\theta), \sin(\theta))$ are the upward unit normal vectors and the unit tangential vectors to the free top surface and $\theta = \arctan (y_{\alpha}/x_{\alpha})$ are the tangent angles on the free top surface. Then $T$ satisfies
\begin{eqnarray}
T_{\alpha} - \theta_{\alpha} U = 0.
\end{eqnarray}
Integrating this relation determines the tangential velocity $T$ in depending on $\theta$ and $U$ up to a constant. Since arc length parametrizations are invariant under translations, this constant can be set to $0$ without loss of generality. This implies
\begin{equation} \label{T} 
T(\alpha,t) = \int_{-\infty}^{\alpha} \theta_{\alpha}(\beta,t)U(\beta,t) \,d\beta\,.
\end{equation}
The normal velocity $U$ is governed by the incompressible Euler's
equations (\ref{euler})--(\ref{incompr}) and the boundary conditions
(\ref{surf})--(\ref{bot}).
For irrotational flows the normal velocity
$U$ can be expressed as a function of the physical tangential velocity $v$ of the fluid particles on the free top surface, whose evolution is determined by (\ref{euler}). Then, using all the above information, one can derive the following 
closed evolutionary system:
\begin{eqnarray}
y_{t} &=& U \cos \theta + {T} {y}_{\alpha}\,, \label{ytilde}\\[2mm]
{v}_{t} & =& - g {y}_{\alpha} + bgh^2 {\theta}_{{\alpha}\, {\alpha}} -  {\delta} {\delta}_{{\alpha}} + U {\theta}_{t}\,, \label{vtilde}\\[2mm]
\theta_{t} &=& U_{\alpha}  - {T} {\theta}_{{\alpha}}\,, \label{thetatilde}\\[2mm]
{\delta}_{{\alpha}\, t} &=& -  {c}{\theta}_{{\alpha}} + bgh^2 {\theta}_{{\alpha}\,{\alpha}\,{\alpha}} -  ({\delta}{\delta}_{{\alpha}})_{{\alpha}} + (U_{\alpha}+ v {\theta}_{{\alpha}})^{2}\,, \label{deltaalphatilde}
\end{eqnarray}
where
\begin{equation}
{c} = {U}_{t} + v {\theta}_{{t}}+ {\delta} {U}_{{\alpha}}+ {\delta} v {\theta}_{{\alpha}} + g \cos{\theta} = -\nabla p \cdot \hat{n}\,, \label{ctilde}
\end{equation}
\begin{equation}
{\delta} = {v} - {T}\,. \label{deltatildevtilde}
\end{equation}
The equations \eqref{thetatilde} and \eqref{deltaalphatilde} are included because they have better regularity properties than the evolution equations for the spatial derivatives of $y$ and $v$. 

The main advantage of this evolutionary system is that in the case with surface tension, i.e., with $b>0$,
the term with the most derivatives in \eqref{ytilde}--\eqref{deltaalphatilde} is linear.

The analog of the arc length formulation in the case of the 3--D water wave problem is the use of isothermal coordinates; see \cite{AM09}.
\bigskip

A further possibility to formulate the 2D--water wave problem is to express the evolution equations in holomorphic coordinates. In doing so, one has to choose a conformal parametri\-za\-tion of $\overline{\Omega}(t)$. Moreover, it is convenient to extend the functions $\eta$ and $\phi$ holomorphically. This formulation was first used by Nalimov \cite{Na74} and Ovsjannikov \cite{O74} and further developed by Dyachenko, Kuznetsov, Spector and Zakharov \cite{DKSZ96}, Wu \cite{Wu97} as well as recently by Harrop--Griffiths, Hunter, Ifrim and Tataru \cite{HIT14, HaIT16}. 

An essential advantage of holomorphic coordinates is that in these coordinates the Dirichlet to Neumann operator is given in terms 
of the Hilbert transform in case of infinite depth of water and in terms of the operator $K_0$ with the symbol 
\begin{equation} \label{K0}
\widehat{K}_0(k)= -i \tanh(hk)
\end{equation}
for $k \in \R$ in case of 
water with finite depth $h$.
\bigskip

Another alternative approach to formulate the water wave problem is to derive an abstract coordinate independent evolutionary system for differential geometric quantities related to the evolving free surface like, for example, its mean curvature. 
Such a system involves differential geometric operators like, for example, the surface Laplace--Beltrami operator and the Riemann curvature tensor. This geometric approach was chosen by Shatah and Zeng \cite{SZ08a, SZ08b, SZ11} and for a related problem also by Beyer and G\"unther \cite{BG98, BG00}.

This approach is based on the fact that the solutions of \eqref{euler}--\eqref{bot} can be interpreted as the geodesic flow with respect to the potential energy, the kinetic energy and in case of surface tension also the surface energy on the infinite dimensional Riemannian manifold of volume--preserving homeomorphisms of $\Omega$. In this variational framework, the boundary conditions \eqref{surf}--\eqref{bot} appear as natural boundary conditions. The interpretation of the incompressible Euler's equation as the geodesic equation on the infinite dimensional manifold of volume--preserving diffeomorphisms was first given by Arnold \cite{A66}. 

The main advantage of a coordinate independent formulation of the water wave 
problem is that it allows one to switch between different coordinate systems. 

\subsection{The linearized water wave problem}

Both for the well--posedness theory and the qualitative theory of the water wave problem the properties of the linearized problem play an important role.
Linearizing \eqref{surf2}--\eqref{pot} around the trivial solution $\eta=v_1=0$
yields the system
\begin{align} 
\tilde{\eta}_t &=  {K}_0 \tilde{v}_1 \,,   \label{surf2lin} 
\\[2mm] 
(\tilde{v}_1)_t &=  - g \tilde{\eta}_x + bgh^2 \tilde{\eta}_{xxx} \,, \label{potlin}
\end{align}
where the operator $K_0$ is defined by \eqref{K0}. The system is hyperbolic but for $b=0$ not strictly hyperbolic.

The system  \eqref{surf2lin}--\eqref{potlin} has plane waves solutions of the form
\begin{equation} \label{linsol1}
\Big(\!\begin{array}{c} \tilde{\eta} \\ \tilde{v}_1 \end{array}\!\Big)  (x,t)=
\Big(\!\begin{array}{c} A \\ K_0^{-1}A_t \end{array}\!\Big) (x,t)\,, 
\end{equation}
with
\begin{equation} \label{linsol2}
A(x,t) = A_0 \sin(k(x + c(k) t)+x_0)\,, 
\end{equation}
where $A_0,k,x_0 \in \R$ are fixed and 
\begin{equation} \label{disprel}
kc(k) = \omega(k) = \pm \sqrt{g(|k| + bh^2|k|^{3})\tanh(h|k|)}\,.
\end{equation}
Hence, the system \eqref{surf2lin}--\eqref{potlin} is dispersive, that means that waves with different wave numbers $k$ propagate with different 
velocities $c(k)$ determined by the dispersion relation $\omega$ which is given by \eqref{disprel}.  

With the help of the 
Fourier transform the system \eqref{surf2lin}--\eqref{potlin} can be uniquely solved for all $t \in \R$ if $\tilde{\eta}|_{t=0}$ and $\tilde{v}_1|_{t=0}$ are regular enough that the Fourier transform can be applied.

The system \eqref{surf2lin}--\eqref{potlin} possesses the conserved quantity
\begin{equation} \label{Elin}
\mathcal{E}_{lin} = \frac12 \int_{\R}  \tilde{\eta} K_0^{-1} g (- \partial_x + bh^2 \partial_x^3) \tilde{\eta} \, dx + \frac12 \int_{\R} \tilde{v}_1^{2} \, dx\,.
\end{equation} 
%where $\widehat{K}_0^{-1}(k) = (\widehat{K_0}(k))^{-1}$. 
$\sqrt{\mathcal{E}_{lin}}$ is 
equivalent to the Sobolev norm $\|\tilde{\eta}\|_{H^{j}} + \|\tilde{v}_1\|_{H^{0}}$,  with $j=1/2$ if $b=0$, and $j=3/2$ if $b > 0$.

By replacing $K_0$ in \eqref{surf2lin}--\eqref{potlin} by the Hilbert transform one gets the linearized system for the case of infinite depth of water. Then the resulting dispersion
relation is 
\begin{equation} \label{disprel2}
\omega(k) = \pm \sqrt{g(|k| + bh^2|k|^{3})}\,.
\end{equation} 
Its second derivative $\omega''$ has no zeros if and only if either $b=0$ or $g=0$. Hence,
an application of the method of stationary phase yields that any sufficiently regular solution of the linearized 2--D water wave problem with infinite depth and either no surface tension or no gravity satisfies the $L^{1}-L^{\infty}$ dispersive decay estimate 
\begin{align} \label{lindisdec}
\|\tilde{\eta}(t)\|_{L^{\infty}} &\leq C (1+t)^{-1/2} (\|\tilde{\eta}(0)\|_{W^{2,1}} + \|\tilde{\eta}_t(0)\|_{W^{2,1}})
\end{align}
for all $t \geq 0$, where $C>0$ is a constant depending neither on $\tilde{\eta}$ nor on $t$. The initial data have to be bounded in the Sobolev space $W^{2,1}$ and not only in $L^{1}$ due to the fact that $\omega''(k) \to 0$ for $ |k| \to \infty$. For the linearized 2--D water wave problem with infinite depth, surface tension and gravity, the $L^{1}-L^{\infty}$ decay rate is only of order $\mathcal{O}((1+t)^{-1/3})$.

The dispersive properties of the linearized 2--D water wave problem with finite depth mainly differ from the case of infinite depth for small wave numbers $k$ since for $|k| \to 0$ the first derivative of the dispersion relation \eqref{disprel} remains bounded and the second derivative tends to $0$.

The dispersion relations of the linearized 3--D water wave problem are also of the form \eqref{disprel} and \eqref{disprel2}, with $k \in \R^2$. Due to the higher space dimension, the $L^{1}-L^{\infty}$ dispersive decay is stronger and has a decay rate of order $\mathcal{O}((1+t)^{-\alpha})$ with $\alpha > 1/2$. For example, the $L^{1}-L^{\infty}$ decay rate of the linearized 3--D water wave problem with infinite depth, gravity and no surface tension is of order $\mathcal{O}((1+t)^{-1})$.

\section{Well--posedness results}  \label{sec3}

\subsection{Local well--posedness}
 
There is a vast literature on the local well--posedness of the water
wave problem. In the case of the 2--D water wave problem with finite depth of water in Eulerian coordinates, local well--posedness means the following. Given the initial data $(\eta_0, (v_1)_0)$ belonging to a physically reasonable subset of a Banach space $X$. Then there exists a $T>0$
and a unique solution $(\eta, v_1) \in C([0,T), X)$ of (\ref{surf2})--(\ref{pot}) with $(\eta|_{t=0}, v_1|_{t=0})=(\eta_0, (v_1)_0)$. Moreover, the solution $(\eta, v_1)$ depends continuously on the initial data. In the other cases, the notion of local well--posedness is analogous.

For analytical initial data, local well--posedness
of the 2D--water wave problem with finite depth of water and without surface tension
was proven in holomorphic coordinates by Ovsjannikov \cite{O74}, in Lagrangian coordinates by Shinbrot \cite{Sh76} and in  
Eulerian coordinates by Kano and Nishida \cite{KN79}. 

Local well--posedness
results for initial data in Sobolev spaces were first obtained in Lagrangian coordinates.
In the case of infinite depth of water, Nalimov \cite{Na74}
proved local well--posedness of the 2--D water wave problem without surface tension for sufficiently small initial data
in Sobolev spaces. Later, Wu \cite{Wu97, Wu99} established local well--posedness of the 2--D and the  3--D water wave problem without surface tension for general initial data in Sobolev spaces. In the case of finite depth of water, 
Yosihara \cite{Y82, Y83} proved local well--posedness of the 2--D water wave problem with and without surface tension for
sufficiently small initial data in Sobolev spaces. Craig \cite {Cr85} showed that the solutions of the 2--D water wave problem without surface tension with sufficiently small initial data in Sobolev
spaces exist even on a longer time interval. This
result was extended by Schneider and Wayne \cite{SW00, SW02a} for the 2--D water wave problem with and without surface tension.

In Eulerian coordinates, local well--posedness of the water wave problem with finite depth and without surface tension for initial data in Sobolev spaces was proven in two dimensions 
by Lannes \cite{L05} and in three dimensions by Alvarez--Samaniego and Lannes \cite{AL08}, where in the latter article a longer interval of existence for sufficiently small initial data was obtained. In \cite{L13}, Lannes gave a simplified and refined proof of these local well--posedness results, which was influenced by Iguchi's proof of the local well--posedness of the 2--D water wave problem with finite depth and with surface tension for initial data in Sobolev spaces 
\cite{I01, I07a}, where in the latter article also a longer interval of existence for sufficiently small initial data was obtained. Analogous results for the 3--D water wave problem with finite depth and surface tension were proven by Ming, Zhang and Zhang \cite{MZ09, MZZ12}. Moreover, local well--posedness results and low regularity theory for the water wave problem in arbitrary dimension with finite or infinite depth and both with and without surface tension were provided by Alazard, Burq and Zuily \cite{ABZ11, ABZ14, ABZ16}.

In the arc--length formulation and in isothermal coordinates, respectively, Ambrose and Masmoudi 
\cite{AM05, AM09} proved local well--posedness of the 2--D and the 3--D water wave problem with and without surface tension in Sobolev spaces. They presented the details of the proof for the case
of infinite depth, but, as they mentioned in \cite{AM05}, the proof works analogously in the case of finite depth. 

In holomorphic coordinates, local well--posedness of the 2--D water wave problem in Sobolev spaces was established by Harrop--Griffiths, Hunter, Ifrim and Tataru in the cases of infinite depth and no surface tension \cite{HIT14},
infinite depth and surface tension but no gravity \cite{IT14b} as well as
finite depth and no surface tension \cite{HaIT16}.

In the abstract geometric formulation, local well--posedness of the water wave problem in Sobolev spaces in arbitrary dimension for finite or infinite depth
and both with and without surface tension was proven by Shatah and Zeng \cite{SZ11}.

Local well--posedness of the water wave problem as a vanishing viscosity limit of the free boundary problem for the Navier--Stokes equations was shown for finite depth and surface tension by Schweizer \cite{Schw05} and in the case without surface tension by Masmoudi and Rousset \cite{MR12}.

For the water wave problem with vorticity, the most general local well--posedness results obtained so far were proven by Christodoulou and Lindblad \cite{CL00}, Coutand and Shkoller \cite{CS07}, Ogawa \cite{O09}, Shatah and Zeng \cite{SZ11} and Castro and Lannes \cite{CL15}.
\bigskip

The general strategy in the proofs of the above local well--posedness results for the water wave problem is to reformulate the water wave problem in a more convenient way, for example, as a quasilinear evolutionary problem, and to derive appropriate energy estimates such that standard techniques for establishing local well--posedness like fixed point iterations or approximating solutions of the exact evolutionary system by solutions of a mollified version of the system can be applied.

To get an idea how an appropriate energy can be constructed it is instructive to consider the formulation \eqref{surf3}--\eqref{pot3} of the 2--D water wave problem. Since the linearization of \eqref{surf3}--\eqref{pot3} around the trivial solution $\eta=v_1=0$ possesses the conserved quantity \eqref{Elin} and since the transport terms in \eqref{surf3}--\eqref{pot3} can be controlled in the Sobolev space $H^{m}$, $m > 3/2$, by using 
$\langle \partial_x^{m} (v_1 \partial_x f), \partial_x^{m} f \rangle_{L^2} \leq C \|v_1\|_{H^{m}} \|f\|_{H^{m}}^2$, with a constant $C>0$ depending neither on $v_1$ nor on $f$, it is natural that an appropriate energy $\mathcal{E}$ for the 2--D water wave problem contains
terms which are related to
\begin{align}
E_s =&\; \frac12 \sum_{\ell=0}^s \Big( \int_{\R} a\, \partial_x^{\ell} \eta K_0^{-1}(-\partial_x) \partial_x^{\ell} \eta \, dx +
bgh^2 \int_{\R} \partial_x^{\ell+1} {\eta} K_0^{-1}(- \partial_x)\partial_x^{\ell+1} {\eta} \, dx \nonumber \\[2mm]
& \qquad\quad\;  + \int_{\R} (\partial_x^{\ell} v_1)^2\, dx \Big)\,.  
\end{align}
If $b>0$, then 
the mapping $(\eta,v_1) \mapsto E_s(\eta,v_1)$ is positive definite for all $(\eta,v_1) \in H^{s+3/2} \times H^{s}$ for sufficiently large $s$. For $b=0$ 
the mapping $(\eta,v_1) \mapsto E_s(\eta,v_1)$ is positive definite for all $(\eta,v_1) \in H^{s+1/2} \times H^{s}$ for sufficiently large $s$
if and only if 
$a > 0$, which is equivalent to the validity of the so--called Taylor--sign condition $ \nabla p \cdot \hat{n} <0$. The Taylor--sign condition was discovered in \cite{T50}.

In many proofs of local well--posedness results in Sobolev spaces an energy $\mathcal{E}$ is constructed in a way that it satisfies the differential inequality
\begin{equation}
\frac{d}{dt} \mathcal{E} \leq f(\mathcal{E})\,,
\end{equation}
with a locally Lipschitz continuous function $f$, such that Gronwall's lemma yields that 
$\mathcal{E}$ is bounded at least for a short time span.

It is shown in \cite{Wu97, Wu99, L13} that the Taylor--sign condition holds for the irrotational water wave problem with infinite depth and with finite depth for a certain class of bottoms containing the flat bottom.
If this condition is violated, then the water wave problem is ill--posed for $b=0$, which was shown in \cite{E87}.
This fact can also be understood from the theory of general hyperbolic systems. Since system \eqref{surf3}--\eqref{pot3} is hyperbolic but for $b=0$ not strictly hyperbolic, the validity of $a > 0$ is a necessary condition for the 
well--posedness of the system in the case $b=0$; see, for example, \cite{FS71, Cr87}. 
If the Taylor--sign condition also holds for $b>0$, then the existence time 
of the solutions does not depend on $b$.

We finish the discussion of the local well--posedness results with the remark that the water wave problem is not completely integrable 
since it has only finitely many conservation laws, see \cite{BO82}. Hence, local well--posedness of the water wave problem cannot be shown with the help of an inverse scattering scheme.

\subsection{Global well--posedness}
For initial value problems from fluid dynamics, it is often a highly nontrivial task to find out how long solutions exist. The most famous example is the question of existence and smoothness of solutions to the three--dimensional incompressible Navier--Stokes equations for all positive times, which is one of the seven Millennium Prize Problems.

For the water wave problem without vorticity,  various almost global and global existence results  were obtained for sufficiently small initial data in Sobolev spaces during the last years. Here almost global existence means that the solution exists for a time span of order $\mathcal{O}(e^{c/\epsilon^n})$ with appropriate constants $c,n >0$ if the initial data are of order $\mathcal{O}(\epsilon)$ for sufficiently small $\epsilon >0$ in a suitable norm and global existence means that the solution exists for all positive times. 

In the case of infinite depth of water, there are the following results.
For the 3--D water wave problem, global
existence of solutions 
was proven in Eulerian coordinates by Germain,
Masmoudi and Shatah in the cases without surface tension \cite{GMS09a, GMS09b} and with surface tension but no gravity \cite{GMS15} as well as by Deng, Ionescu, Pausader and Pusateri \cite{DIPP16} in the case with surface tension and gravity. 
In Lagrangian coordinates, global existence of solutions to the 3--D water wave problem was shown by Wu \cite{Wu09b} in the case without surface tension.
For the 2--D water wave problem without surface tension, almost
global existence of solutions was shown by Wu \cite{Wu09a} in Lagrangian coordinates and by Hunter, Ifrim and Tataru \cite{HIT14} in holomorphic coordinates. Global existence results for the 2--D water wave problem were established in the case without surface tension by Ionescu and Pusateri \cite{IP14} and by Alazard and Delort \cite{AD13, AD15} in Eulerian coordinates and by Ifrim and Tataru \cite{IT14a} in holomorphic coordinates as well as in the case with surface tension but no gravity by Ionescu and Pusateri \cite{IP14b} in Eulerian coordinates and by Ifrim and Tataru \cite{IT14b} in holomorphic coordinates. In the above global existence results, also asymptotic expansions of the solutions for $t \to \infty$ were provided.

For finite depth of water, global existence of solutions to the 3--D water wave problem without vorticity having sufficiently small initial data in Sobolev spaces was recently proven in Eulerian coordinates by Wang in the cases without surface tension \cite{Wa15a, Wa15b} and with surface tension but no gravity \cite{Wa16}.
\bigskip

The strategy in all these almost global and global existence proofs 
is a generalization of a procedure which was developed for proving global well--posedness of several classes of nonlinear wave equations and nonlinear dispersive equations for sufficiently small initial data. This procedure consists in deriving and combining so--called high--order energy estimates and dispersive decay estimates.

The simplest version of a high--order energy estimate is an a priori estimate of the form 
\begin{align} \label{hee}
\|U(t)\|_{W^{m,2}} &\leq C_1 \|U(0)\|_{W^{m,2}}\, \text{exp} \Big(C_1 \int_0^t \|U(s)\|_{W^{r,\infty}}^{\beta}\,\mathrm{ds} \Big) 
\end{align}
for $t \geq 0$, where $U(t)$ is a quantity depending on the respective solution of the considered equation (and sometimes also on some of its derivatives) at time $t$, $C_1>0$ is a constant depending on $m$ but neither on $U$ nor on $t$ and $\beta>0$ is a constant depending on the asymptotic behavior of the nonlinearity of the considered equation as $U \to 0$. Moreover, the Sobolev index $r \geq 0$ is independent of the Sobolev index  $m \geq 0$.  Hence, in order to control $\|U(t)\|_{W^{m,2}}$, one also needs to control $\|U(t)\|_{W^{r,\infty}}$, which shall be done with the help of a dispersive decay estimate.

The simplest version of a dispersive decay estimate is an a priori estimate of the form 
\begin{align} \label{dde}
\|U(t)\|_{W^{r,\infty}} &\leq C_2 (1+t)^{-\alpha} \|U(0)\|_{W^{\ell,1}}
\end{align}
for $t \geq 0$ and sufficiently small $\|U(0)\|_{W^{\ell,1}}$, with $\ell \geq 0$, where $C_2>0$ is a constant depending on $r \geq 0$ but neither on $U$ nor on $t$. Moreover,
$\alpha >0$ is a constant depending on the space dimension, the linear dispersion relation of the considered equation and the asymptotic behavior of the nonlinearity of the considered equation as $U \to 0$.  

Combining the high--order energy estimate \eqref{hee} with the dispersive decay estimate \eqref{dde} 
yields a global a priori bound for $\|U(t)\|_{W^{m,2}}$ if $-\alpha \beta < -1$. This bound allows one to continue a local solution of the considered equation for all positive times. For a more detailed presentation of this procedure and for some examples of equations to which the procedure can be successfully applied, we refer, for example, to \cite{KP83}.

The best dispersive decay rates that can be expected for the water wave problem are of the same order as the dispersive decay rates of the linearized water wave problem. Since the evolution equations of the water wave problem contain quadratic terms, this implies that a combination of an optimal dispersive decay estimate of the form \eqref{dde} with a high--order energy estimate of the form \eqref{hee} would yield in the case of the 2--D water wave problem with infinite depth and either gravity or surface tension an integrand of order $\mathcal{O}((1+s)^{-1/2})$ 
on the right--hand side of \eqref{hee}, which would not be sufficient to derive a global a priori bound for $\|\eta(t)\|_{W^{m,2}}$ and $\|v_1(t)\|_{W^{m,2}}$.
However, if a normal--form transform could be performed to eliminate the quadratic 
terms and to replace them by cubic terms, one would obtain an integrand of order $\mathcal{O}((1+s)^{-1})$ such that almost global existence results and with the help of more subtle techniques also global existence results could be proven. 

In the case of the 3--D water wave problem with infinite depth or with finite depth and either gravity or surface tension the corresponding integrand in \eqref{hee} would be of order $\mathcal{O}((1+s)^{-\gamma})$ with $\gamma >1$ after performing a normal--form transform
such that establishing global existence results would be somewhat easier than in two dimensions.

In realizing the above strategy of deriving and combining high--order energy estimates and dispersive decay estimates for the water wave problem various difficulties arise, in particular, handling occurring resonances, tackling the loss of regularity in the nonlinearity and finding appropriate norms which can take the role of the norms used in \eqref{hee}--\eqref{dde}.  

These difficulties are overcome with the help of a combination of variants of Klainerman's vector fields method \cite{K85} with generalizations of Shatah's normal--form transform method \cite{Sha85} realized in a distinct way 
by the different authors. Wu constructed an implicit change of coordinates which  agrees quadratically with Shatah's normal--form transform and is derived by using complex analysis in the two dimensional case and Clifford analysis in the three dimensional case. Germain, Masmoudi and Shatah developed the so--called space--time resonance method, Alazard and Delort the normal--forms paradifferential method, Hunter, Ifrim and Tataru the modified energy method and the method of testing by wave packets and Ionescu and Pusateri the quasilinear I--method. 
\bigskip

Since the optimal dispersive decay rate of the 2--D water wave problem with infinite depth, gravity and surface tension is only of order $\mathcal{O}((1+s)^{-1/3})$, the above strategy is not sufficient to prove almost global or global existence of solutions, which is still an open problem.

Global or almost global existence of solutions to the 2--D water wave problem with finite depth and to the 3--D water wave problem with finite depth, gravity and surface tension are also open problems.
A reason for this is the fact that dispersive estimates as above cannot be expected because the 2--D water wave problem with finite
depth and the 3--D water wave problem with finite depth, gravity and surface tension possess solitary waves solutions, which are localized traveling waves of permanent form; see below. Answering the question of global existence of solutions in the neighborhood of solitary wave profiles is also of interest in the context of investigating long time stability 
of solitary waves solutions. 

Global existence for arbitrarily large initial data is not true in general since singularities can occur in finite time \cite{CCFGL12, CCFGG12, CCFGG13}. The most
common singularity is of course wave breaking. However, developing a general theory of wave breaking is still an open problem; see \cite[Chapter 6]{C11} for further information.

\section{Reduced models and their validity}

Concerning the qualitative behavior of the solutions, the full
water wave problem -- even in two dimensions -- is extremely complicated to analyze. A qualitative understanding of the solutions to the full water wave problem being usable for practical applications does not seem within reach for the near future, neither analytically nor numerically. Therefore, it is important to approximate the full model by suitable reduced model equations whose solutions have similar but more easily accessible qualitative properties. In the following, we present various reduced models. For simplicity we restrict ourselves again to the 2--D water wave problem. For the corresponding reduced models for the 3--D water wave problem and their properties, we refer, for example, to \cite{SW02b, L13, T14}. 

The reduced models which we will present are used to describe the motion of some typical wave profiles occurring in nature, namely, Airy waves, Stokes waves, cnoidal waves and solitary waves. Each of these wave profiles appears in typical parameter regimes for the ratio of the depth of the water, the wave length and the amplitude of the wave; see, for example, \cite[Chapter 15]{LM76}.
The reduced models can be derived by multiple scaling analysis.

From now on, let space and time in the water wave problem be rescaled in such a
way that $h=1$ and $g=1$.

\subsection{The linear wave equation}
Inserting the multiple scaling ansatz
\[
\Big(\!\begin{array}{c} \eta \\ u_1 \end{array}\!\Big)  (x,t)=
\varepsilon^{\alpha} \Big(\!\begin{array}{c} A \\ B \end{array}\!\Big) 
(\varepsilon x, \varepsilon t) 
+ \mathcal{O} (\varepsilon^{\alpha+1}) 
\]
with $0 < \varepsilon \ll1$, $ A, B: \R^2 \to \R$ and $\alpha >2$,
into \eqref{surf2}--\eqref{pot}, expanding the operator $\mathcal{K}$ with respect to $\varepsilon$ and equating the terms with the lowest power of $\varepsilon$
yields that $A$ has to satisfy in lowest order with respect to $\eps$ the linear wave equation 
\begin{align} \label{waveeq}
A_{\tau\tau} &= A_{\xi\xi} \,,
\displaystyle
\end{align}
where $\tau = \varepsilon t$ and $\xi = \varepsilon x $. 
The linear wave equation is the simplest and most used reduced model for the water wave problem. It is helpful to explain many physically relevant linear wave phenomena.
The linear wave equation \eqref{waveeq} has solutions of the form
\begin{align} 
A(\xi,\tau) &= A_0 \sin(k(\xi \pm \tau)+\xi_0) 
\end{align}
with $A_0,k,\xi_0 \in \R$. These solutions are called Airy waves. They were first derived in the context of water waves by Airy \cite{A41} in 1841.

\subsection{The Korteweg--de Vries equation}
Now, we present nonlinear reduced models for the 2--D water wave problem. 
The most famous nonlinear reduced models are the Korteweg--de Vries (KdV) equation and the Nonlinear Schr\"odinger (NLS) equation.
First, we discuss the KdV equation.
By inserting the ansatz
\[
\begin{pmatrix} \eta \\ v_1 \end{pmatrix} (x,t)=
\varepsilon^2 A \left(\varepsilon (x \pm t), \varepsilon^3t\right) 
\begin{pmatrix} 1 \\ \mp 1  \end{pmatrix} + \mathcal{O} (\varepsilon^3)
\]
with $0 < \varepsilon \ll1$ and $ A: \R^2 \to \R$ into
\eqref{surf2}--\eqref{pot}, expanding the operator $\mathcal{K}$ with respect to $\varepsilon$ and equating the terms with the lowest power
of $\eps$ one obtains that $A$ has to satisfy in lowest order with respect to $\eps$ the KdV equation
\begin{equation} \label{kdv}
A_{\tau} = \pm \Big(\frac16-\frac{b}{2}\Big) A_{\xi\xi\xi} \pm \frac32 A A_{\xi}\,,
\end{equation}
where $\tau = \varepsilon^3 t$ and $\xi = \varepsilon (x\pm t)$.

\vspace*{0.35cm}
\begin{figure}[htbp]
\epsfig{file=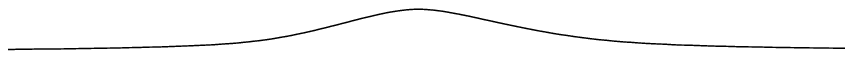,width=15cm,height=1.5cm,angle=0}
\vspace*{-1.8cm}

\hspace*{3.7cm}
\epsfig{file=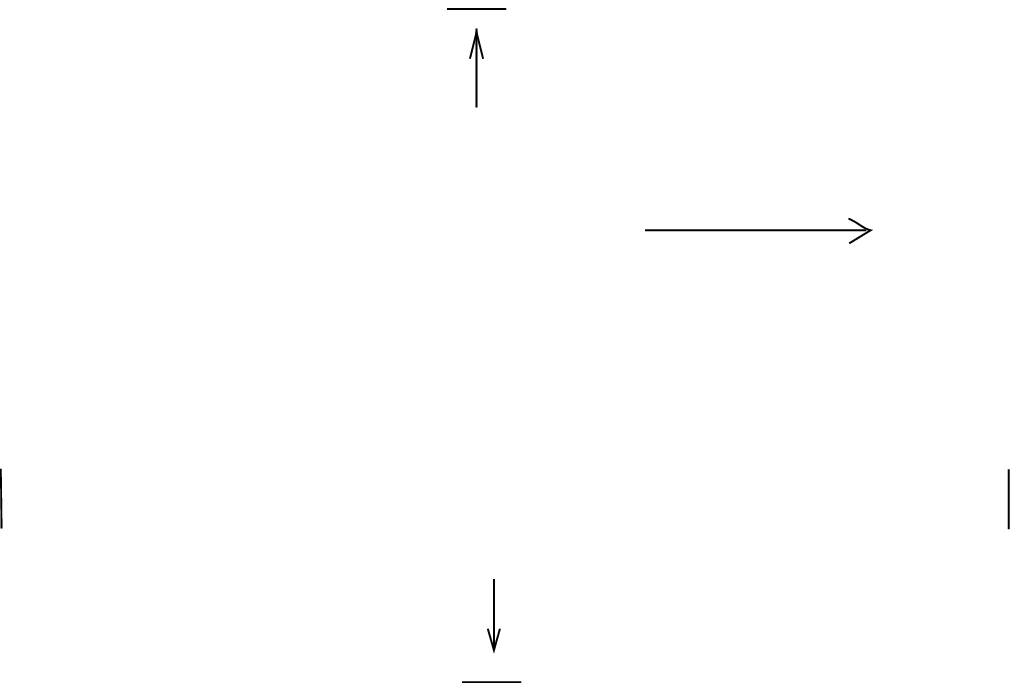,width=7.5cm,height=2.5cm,angle=0}
\vspace*{-5.34cm}
%\hspace*{9.7cm}
%$c_\mathrm{g}$  
\vspace*{3.2cm}

%\hspace*{9.7cm} 
%$c_\mathrm{p}$

\vspace*{0.6cm}
\hspace{2.1cm}\hspace{5.0cm}
$\varepsilon^2$
\vspace*{0.1cm}

\hspace{4.5cm}
$1/\varepsilon$%\hspace{4.9cm}$1/\varepsilon$
\hspace{2.7cm}
\vspace*{0.75cm}

\caption{{\small A solution $\eta$ to \eqref{surf2}--\eqref{pot} of order $\mathcal{O}(\varepsilon^2)$ moving to the right (or to the left) whose characteristic length scale is of order $\mathcal{O}(1/\varepsilon)$
is approximately described by a solution $A$ to the KdV equation \eqref{kdv}.} \label{fig2} }
\end{figure}

A similar formal perturbation analysis was first done by Boussinesq \cite{Bou71}, \cite {Bou77}, Rayleigh \cite{Ray76} and Korteweg and de Vries \cite{KdV95} in the second half of the nineteenth century.

For $b = \frac13 + 2\nu \epsilon^2$ one gets, by 
making the ansatz
\[
\begin{pmatrix} \eta \\ v_1 \end{pmatrix} (x,t)=
\varepsilon^4 A \left(\varepsilon (x \pm t), \varepsilon^5t\right) 
\begin{pmatrix} 1 \\ \mp 1  \end{pmatrix} + \mathcal{O} (\varepsilon^5)
\]
with $0 < \varepsilon \ll1$ and $ A: \R^2 \to \R$ the Kawahara equation
\begin{equation} \label{kara}
\partial_{\tau} A = \mp \nu \partial_{\xi}^3 A \pm \frac{1}{90} \partial^5_{\xi} A \pm \frac32 A \partial_{\xi}A\,,
\end{equation}
where $\tau = \varepsilon^5 t$ and $\xi = \varepsilon (x \pm t)$
as an approximation equation; see, for example, \cite{SW02a}.

The KdV equation and the Kawahara equation can be derived as
approximation equations only in the case of finite depth of water since the expansion of the operator $\mathcal{K}$ with respect to $\varepsilon$ is different for infinite depth of water.

The important point about the KdV equation is that the effects induced by the nonlinear term
are in some equilibrium with the dispersive effects, and so the KdV equation
has a family of solitary waves, which are spatially localized traveling waves of permanent form, as
solutions; see \cite{KdV95}. For $b >\frac13$, shape and propagation
behavior of the solitary wave solutions to \eqref{kdv} are different
than in the regime $b < \frac13$; see, for example, \cite{SW02a}. Furthermore, the
KdV equation possesses N--solitons as solutions; see \cite
{ZK65}. N--solitons are N waves of permanent form that go through each
other without changing their shape asymptotically. The Kawahara
equation possesses solitary wave solutions at least if $\nu=1$; see \cite{YT81, Il90}.

The above perturbation analysis suggests that the soliton dynamics of
the KdV equation is at least approximately present in the 2--D water wave problem. Indeed, solitary waves were observed experimentally in a canal between Edinburgh and Glasgow by Russell in 1834; see \cite{R44}. 

Moreover, the KdV equation also has periodic solutions, which are called cnoidal waves. Cnoidal waves have sharper crests and a flatter troughs than Airy waves.

The full water wave problem also possesses exact solitary wave solutions and exact periodic wave solutions. The latter are called Stokes waves and were first derived by Stokes \cite{Sto47} in 1847.  For an overview of the vast literature about existence 
results for these types of waves in different parameter regimes and characteristic properties of their profiles as well as the current state of the art concerning their stability properties, we refer, for example, to the surveys \cite{DI03, G04, G07, CW07, Str10, Sun14}, the book \cite[Chapter 3--5]{C11}
and the recent articles \cite{VW14, MRT15, PS16, GSW16, BGW16}. Furthermore, for results about the form of particle paths beneath Stokes waves and solitary waves,
we refer, for example, to \cite{C06, CE07}, \cite[Chapter 4--5]{C11} and for results about monotonicity properties of the pressure 
beneath Stokes waves and solitary waves to \cite{CS10, CEH11}, \cite[Chapter 4--5]{C11}. 

Many of the results about solitary waves and Stokes waves in the literature were proven with the help of an interplay between subtle methods from diverse areas of analysis, in particular elliptic theory, nonlinear functional analysis, dynamical systems and calculus of variations.

Both Airy waves and solitary waves have also been studied in the context of modeling tsunamis.
For a discussion of the possibilities and limitations of using these waves to describe the dynamics of tsunamis; see, for example, \cite[Chapter 7]{C11} and the references therein.
\bigskip

A rigorous theoretical confirmation that the
soliton dynamics of KdV equation 
and the dynamics of the Kawahara equation are at least approximately
present in the 2--D water wave problem can be given 
by proving error estimates, i.e., estimates of the difference between exact solutions of the 2--D water wave
problem and their approximations obtained via the KdV or the Kawahara
equation. In the case of zero
surface tension, the first approximation proofs were done by Craig \cite{Cr85} in Lagrangian
coordinates and by Kano and Nishida \cite{KN86} in Eulerian coordinates.
However, the approximation result of Kano and Nishida is restricted to
analytical initial data, and, as mentioned in \cite{SW00},  the validity of
the approximation is proven  only for a time interval of order $\mathcal{O}(\eps^{-1})$ in
the original system. This time span is too short to see the physically relevant dynamics of
the KdV equation, for example, the interaction of solitary waves, in the 2--D water wave problem since, as we can see
from \eqref{kdv}, the
characteristic time scale of the KdV approximation is of order
$\mathcal{O}(\eps^{-3})$. 
The result of Craig admits a least a certain class of initial data
in Sobolev spaces and is valid on a time interval of order $\mathcal{O}(\eps^{-3})$ in
the original system. 

Craig's result was improved by Schneider and Wayne
\cite{SW00} in the following sense. They included Sobolev initial data
that are in a neighborhood of solitary wave or multisoliton profiles.
Furthermore, they extended the time span of validity to a longer time
interval of order $\mathcal{O}(\eps^{-3})$ which includes the interaction time
of solitons. Finally, they proved that solutions of the 2--D water
wave problem with small, long--wave, sufficiently smooth and
sufficiently localized
initial profiles split up into two wave packets, one moving to the
left and one moving to the right, where each of these wave packets can
be well approximated by solutions of two decoupled KdV equations.   
In \cite{SW02a} they also proved analogous justification theorems for
the KdV and the Kawahara approximation of the 2--D water wave problem
with surface tension.  
Wright \cite{W05} extended Schneider's and
Wayne's result by adding higher--order correction
terms to the KdV approximation.

Error estimates for the KdV or the Kawahara
approximation yielding slightly improved approximation rates were proven in Eulerian coordinates by Bona, Colin and Lannes \cite{BCL05} in the case without surface tension and by
Iguchi \cite{I07a, I07b} in the case with surface
tension. These error estimates 
are valid for the same initial data and the same time intervals as
Schneider's and Wayne's estimates. Similar approximation results for 
slightly varying bottoms were provided by Chazel \cite{Ch09} in the case without surface tension and by Iguchi \cite{I06, I07a} in the case with surface
tension.

In \cite{D12}, a simplified approximation proof was given, which was performed in the
arc length formulation of the 2--D water wave problem. Due to the structure of the arc length formulation the error
estimates in  \cite{D12} are the only ones being
uniform with respect to the strength of the surface tension, as $b$
and $\eps$ go to $0$. 
%Therefore, the rigorous justification of the KdV approximation
%can be given for the cases with and without surface tension together
%by one proof. 
The proof of the Kawahara approximation goes
analogously. The error estimates 
are valid for the same initial data and the same time intervals as
Schneider's and Wayne's estimates and yield the same approximation rates
as the estimates by Bona, Colin and Lannes and Iguchi, respectively. 
Transferring the error estimates into Eulerian coordinates preserves
their accuracy since in the scaling regimes of the
KdV or the Kawahara approximation the coordinates of the free surface
in arc length parametrization
are very close to its Eulerian coordinates.

Expressed in Eulerian coordinates the approximation result for the KdV equation in \cite{D12} reads as follows.

%\pagebreak

\begin{theorem} \cite[Theorem 1.1]{D12} \label {th11} 
For all $b_0, C_0, \tau_0 > 0$ there exist an
$\epsilon_0 > 0$ such that for all $\epsilon \in \R$ with
$0 < \eps \leq \epsilon_0$ and all $b \in \R \setminus \{\frac13\}$ with  $0
\leq b \leq b_0$  the following is true. Let 
\[
\eta|_{t=0}(x) = \epsilon^2 \Phi_1(\epsilon x), \qquad v_1|_{t=0}(x) = \epsilon^2 \Phi_2(\epsilon x)
\]
with $\|(\Phi_1,\Phi_2)\|_{H_{\xi}^{s+8} \,  \cap
  H_{\xi}^{s+3}(k)} \leq C_0 \eps^{l}$, where  $\xi=\eps x$,\, $s \geq 7$,\,  $k >1$ and $l \geq 0$. Split the initial conditions into 
\[
A_1|_{\tau=0} = \frac12 (\Phi_1 + \Phi_2), \qquad A_2|_{\tau=0} = \frac12 (\Phi_1 - \Phi_2)
\]
and let the amplitudes $A_1= A_1(\xi,\tau)$ and $A_2=A_2(\xi,\tau)$ satisfy 
\[
(A_1)_{\tau} = \Big(\frac{b}{2}-\frac16\Big) (A_1)_{\xi\xi\xi} -\frac32 A_1 (A_1)_{\xi}, \qquad 
(A_2)_{\tau} = \Big(\frac16-\frac{b}{2}\Big) (A_2)_{\xi\xi\xi} +\frac32 A_2 (A_2)_{\xi}.
\]
Then there exists a unique solution of the 2--D water wave problem (\ref{surf2})--(\ref{pot}) with the above initial conditions satisfying
\[
\sup\limits_{t \in [0, {\tau}_0 /\epsilon^3]} 
\left\| 
  \begin{pmatrix}
\eta \\ v_1
\end{pmatrix} \!
(\cdot, t) - 
 \psi(\cdot , t) 
\right\|_{H^{s}_{\xi}  \times H^{s-1/2}_{\xi}} \lesssim  \epsilon^{4+l}
\] 
where 
\[
\psi(x, t) = \epsilon^2 A_1 \left(\epsilon (x-t), \epsilon^3t\right) 
\begin{pmatrix} 1 \\ 1  \end{pmatrix} + \epsilon^2 A_2 \left(\epsilon (x+t), \epsilon^3t\right) 
\begin{pmatrix} 1 \\ -1  \end{pmatrix}.
\]
\end{theorem}

Here $H_{\xi}^s$ denotes the inhomogeneous
  Sobolev space of order $s$ on $\R$ with integration variable ${\xi}$ and  $H_{\xi}^s(k)$  the inhomogeneous weighted
  Sobolev space of order $s$ on $\R$ with integration variable ${\xi}$
  and weight  $\rho^k({\xi})=(1+{\xi}^2)^{k/2}$, i.e., 
$\|u\|_{H_{\xi}^s(k)} = \|u \rho^k\|_{H_{\xi}^s}$. Moreover, $A
\lesssim B$ means $A \leq C B$ for a constant $C>0$.

An error of order $\mathcal{O}(\epsilon^{4+l})$ is small compared with the exact solution 
and its approximation which are both of order $\mathcal{O}(\epsilon^2)$. We note that this fact should not be taken for granted. There are approximation equations which although derived by reasonable formal arguments do not reflect the true dynamics of the corresponding original equations; see, for example, \cite{GS01, SSZ15}. 

The main difficulty in proving the above approximation result is to control the error $R$ made by the approximation $\psi$ on the long time span of order $\mathcal{O}(\epsilon^{-3})$. This problem is handled by constructing an energy $\mathcal{E}= \mathcal{E}(R)$ being equivalent to the square of a suitable Sobolev norm of $R$ for sufficiently small $\varepsilon$ and $R$ in a way that 
\begin{equation} \label{Ein}
\frac{d}{d{t}} \mathcal{E} \lesssim \varepsilon^{3} \Big(\mathcal{E} +1 + \frac{1}{\varepsilon^{2}(1+(\varepsilon t)^{2})^{k/2}} \Big)
\end{equation}
as long as $\mathcal{E} = \mathcal{O}(1)$, where the last summand on the right--hand side reflects the fact that the interaction time of the two localized waves $A_1$ and $A_2$ is only of order $\mathcal{O}(\varepsilon^{-1})$. By using 
\eqref{Ein} and Gronwall's lemma, one obtains the $\mathcal{O}(1)$--bounded\-ness of $\mathcal{E}$ and consequently the desired estimate of $R$ for all $t\in[0,\tau_0/\varepsilon^{3}]$. 
The construction of this energy is governed by the structural properties and the scaling behavior of the linear and the quadratic components of the evolutionary
system \eqref{ytilde}--\eqref{deltatildevtilde} of the water wave problem in the arc length formulation. 

The approximation theorem for the Kawahara equation is similar. 

\subsection{The Whitham system}
The next nonlinear reduced model we present is the Whitham system. 
By inserting the ansatz
\[
\Big(\!\begin{array}{c} \eta \\ u_1 \end{array}\!\Big)  (x,t)=
\Big(\!\begin{array}{c} A \\ B \end{array}\!\Big) 
(\varepsilon x, \varepsilon t)  + \mathcal{O} (\varepsilon) 
\]
with  $0 < \varepsilon \ll1$ and $ A, B: \R^2 \to \R$ into
\eqref{surf2}--\eqref{pot}, expanding the operator $\mathcal{K}$ with respect to $\varepsilon$ and equating the terms with the lowest power
of $\varepsilon$, one obtains a nonlinear system of conservation laws of the form
\begin{align}
\label{wh1}
\partial_T A &= -\partial_X ((1+A) B) \,, \\
\label{wh2}
\partial_T B &= -\partial_X A
-\partial_X \Big(\frac{B^2}{2}\Big) \,,
\end{align}
where $ X = \varepsilon x $,  $T = \varepsilon t $.
This system is called Nonlinear Shallow Water system or Whitham system. 
It was first derived by Saint-Venant 
\cite{SV71a,SV71b} in 1871.

The validity of this reduced model was rigorously justified for the 2--D water wave problem without surface tension by Ovsjannikov \cite{O74, O76} for analytic initial data and by Iguchi \cite{I09} for initial data in Sobolev spaces.

A related system was derived by Whitham \cite{Wh65a, Wh65b, Wh74} in 1965
for the approximate description of slow modulations in time and space of  periodic wave--trains in dispersive wave systems. This system is also called Whitham system.

Moreover, Benjamin and Feir \cite{B67, BF67} as well as Whitham \cite{Wh67} discovered that such wave--trains can become unstable if their wave numbers are sufficiently large. This instability is called 
the Benjamin--Feir instability. The occurrence of the Benjamin--Feir instability was rigorously shown by Bridges and Mielke \cite{BM95}.

In \cite{BDS16}, the Whitham approximation was justified for a model problem with periodic coefficients. The justification of the Whitham system for the approximate description of slow modulations in time and space of periodic wave--trains in the water wave problem is 
still an open problem and subject of current research.

\subsection{The Nonlinear Schr\"odinger equation}
Finally, we discuss the NLS equation. The KdV equation, the Whitham system 
and the NLS equation are the three generic modulation equations for dispersive systems like the water wave problem.

\begin{figure}[htbp]
\epsfig{file=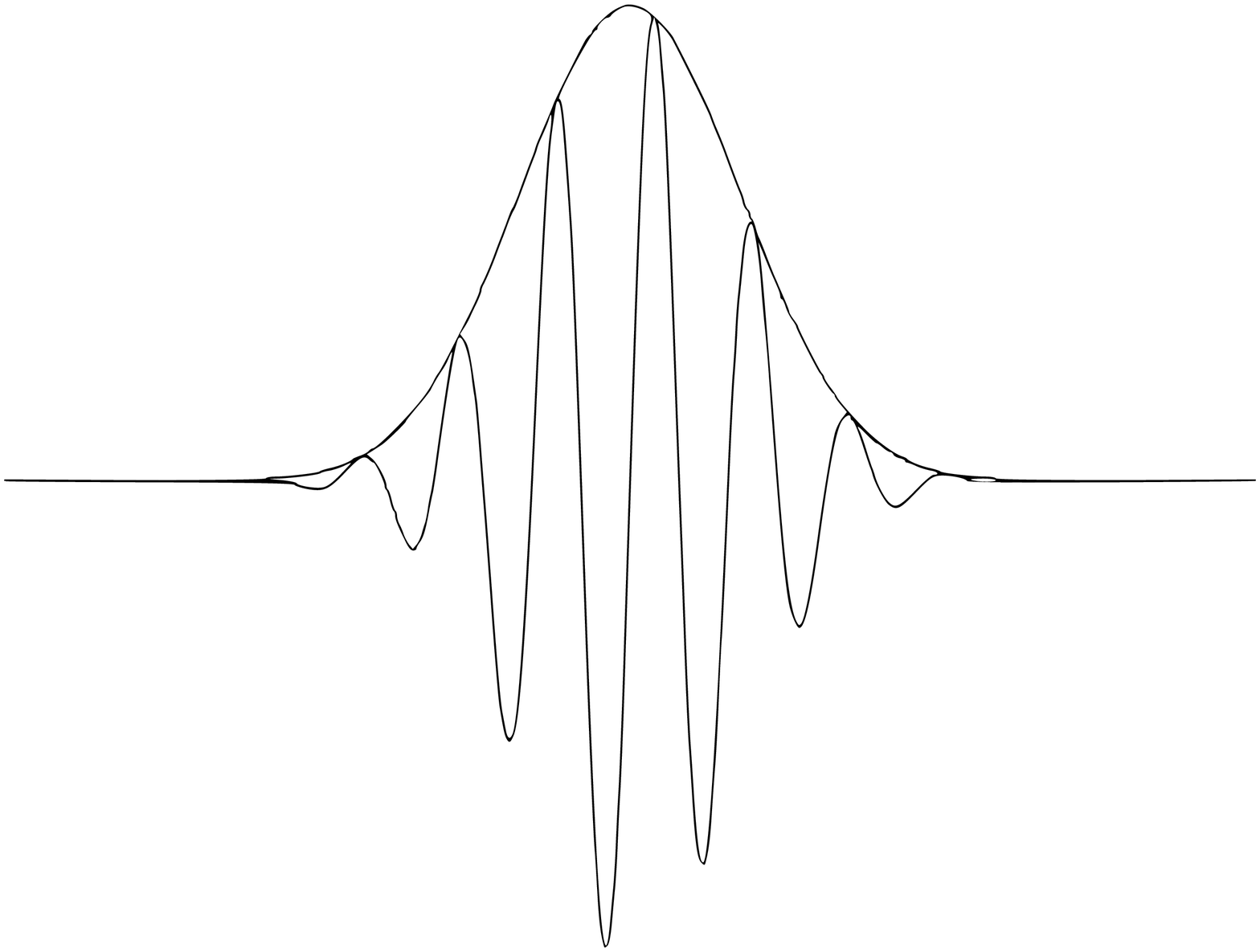,width=15cm,height=6.2cm,angle=0}
\vspace*{-6.3cm}

\hspace{1.9cm}
\epsfig{file=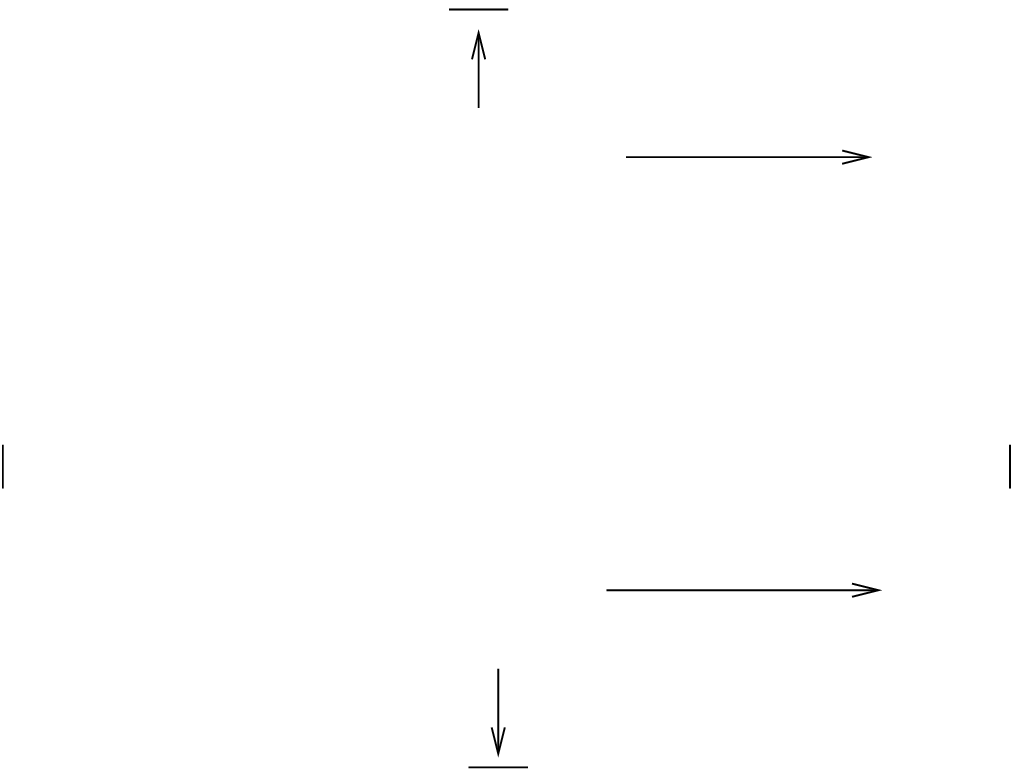,width=10.7cm,height=6.5cm,angle=0}

\vspace*{-5.34cm}
\hspace*{9.7cm}$c_\mathrm{g}$  \vspace*{3.2cm}

\hspace*{9.7cm} $c_\mathrm{p}$

\vspace*{-2.5cm}
\hspace{2.2cm}\hspace{5.0cm}$\varepsilon$
\vspace*{0.3cm}

\hspace{2.65cm}$1/\varepsilon$%\hspace{4.9cm}$1/\varepsilon$
\hspace{2.7cm}
\vspace*{2.6cm}

\caption{{\small An envelope (advancing with the group velocity 
$ c_g $) with characteristic length scale of order $\mathcal{O}(1/\varepsilon)$ of an oscillating wave packet $\eta$ of order $\mathcal{O}(\varepsilon)$ 
(advancing with the phase velocity 
$ c_p= \omega_0/k_0  $) is approximately described by the 
amplitude $ A $ which solves  the NLS equation \eqref{NLS}.} \label{fig3} }
\end{figure}

The ansatz for the NLS approximation is 
\[
\begin{pmatrix} \eta \\ v_1 \end{pmatrix} (x,t)= 
\varepsilon A \left(\varepsilon (x -c_gt), \varepsilon^2t\right) e^{i(k_0x-\omega_0 t)}
\varphi(k_0) + \mathcal{O} (\varepsilon^2) + \mathrm{c.c.}\,,
\]
where $0 < \varepsilon \ll 1$.
Here $ \omega_0 >0$ is the basic temporal 
wave number associated via the linear dispersion relation \eqref{disprel} to the basic spatial wave number $ k_0 > 0$ of the underlying carrier wave $ e^{i(k_0 x - \omega_0 t)}$, that means that $\omega_0= \omega(k_0)$, where the branch of solutions $\omega(k) = {\rm sign}(k) ((k+bk^3)\tanh(k))^{1/2}$ is chosen. Moreover,
$c_g$ is the group velocity, i.e., $c_g=  \omega'(k_0)$, $A$ the complex--valued amplitude, $\varphi(k_0) \in \C^2$ and c.c. the complex conjugate.  This ansatz leads to waves moving to the right; to obtain waves moving
to the left, $\omega_0$ and $c_g$ have to be replaced by  $-\omega_0$ and $-c_g$.
By inserting this ansatz into (\ref{surf2})--(\ref{pot}), one obtains that for an explicitly computable vector $\varphi(k_0)$ the amplitude $A$ has to satisfy at leading order in $\eps$ the NLS equation
\begin{equation} \label{NLS}
A_{\tau} = i \frac{\omega''(k_0)}{2} A_{\xi\xi} +i \nu(k_0) A |A|^2 \,,
\end{equation}
where $\tau = \varepsilon^2 t, \xi = \varepsilon (x -c_g t)$ and $\nu(k_0) \in \R$.
Hence, the NLS equation \eqref{NLS} approximately describes the dynamics of spatially and temporarily oscillating wave packet--like solutions to the 2--D water wave problem; see Figure \ref{fig3}.

\vspace*{0.35cm}

The NLS approximation was first derived by Zakharov \cite{Za68} in 1968. In numerical simulations, the simulation of the evolution of the envelope with the help of the NLS approximation yields a significant reduction of complexity and consequently an increase of efficiency compared to the simulation of the whole wave packet. 
The NLS approximation as well as the Whitham approximation are also of interest in the context of modeling monster waves; see \cite{JO07}.

The NLS equation \eqref{NLS} possesses time--periodic solitary wave solutions if $
\omega''(k_0)$ and $\nu(k_0)$ have the same sign; for a discussion of the values 
of $\nu(k_0)$; see \cite[Chapter 4]{AS81}.

The NLS equation and the KdV equation belong to wider classes of reduced models for the water wave problem; see \cite{L13} for an overview of their derivation and justification results. In these classes of reduced models, the KdV equation and the NLS equation have the advantage of being independent of the parameter $\varepsilon$. Moreover, the one dimensional KdV equation and NLS equation are completely integrable Hamiltonian systems which 
can be explicitly solved with the help of inverse scattering schemes.
\bigskip

In \cite{DSW16}, the NLS approximation was justified for the 2--D water wave problem with finite depth and no surface tension by the following approximation result.

\begin{theorem}
\label{th2} \cite[Theorem 1.1]{DSW16}
Let $b=0$ and $s \geq 7$.  Then for all $k_0 > 0$ 
and for all $C_1,{\tau}_0 > 0$ there exist $ {\tau}_1 > 0 $ and $\epsilon_0 > 0$ 
such that for all solutions $A \in
C([0,{\tau}_0],H^{s}(\R,\C))$ of the NLS equation (\ref{NLS})
with 
$$ 
\sup\limits_{\tau \in [0,{\tau}_0]} \| A(\cdot,\tau) \|_{H^{s}(\R,\C)} \leq C_1
$$  
the following holds.
For all $\epsilon \in (0,\epsilon_0)$
there exists a solution 
$$
(\eta,v_1)\in
C([0,{\tau}_1/\epsilon^2],(H^{s}(\R,\R))^2)
$$
of the 2-D water wave problem
(\ref{surf2})--(\ref{pot}) which satisfies
  $$\sup\limits_{t \in [0,{\tau}_1/\epsilon^2]} \left\| \! \begin{pmatrix}
\eta \\ v_1
\end{pmatrix} \!(\cdot,t) -
    \psi(\cdot,t)\right\|_{(H^{s}(\R,\R))^2} 
\lesssim \epsilon^{3/2}$$
where
$$
\psi(x, t) = \varepsilon A \left(\varepsilon (x -c_gt), \varepsilon^2t\right) e^{i(k_0x-\omega(k_0)t)}
\varphi(k_0) + c.c.\,.
$$
\end{theorem}

The theorem is proven by using Lagrangian coordinates.
The main difficulty of the proof is again to control the error $R$ made by the approximation $\psi$ on a long time span, here of order $\mathcal{O}(\epsilon^{-2})$. 
 This problem is handled by eliminating all terms of order $\mathcal{O}(\epsilon)$ in the evolutionary system of the error by a normal--form transform 
of the form
\begin{equation} \label{inft} 
\tilde{R} := R + \eps
N({\psi},R)\,,
\end{equation}
where $N$ is an appropriately constructed bilinear mapping; compare \cite{Sha85, Kal88}. 
Here, such a normal--form transform can be constructed although the evolution equations of the error possess resonances which cause a small divisor problem. Even though the resulting normal--form transform loses regularity the structure of the error equations allows  one to invert the normal--form transform and to express the error equations
in terms of the new variable $\tilde{R}$.
Due to the properties of the Lagrangian formulation of the water wave problem 
the right--hand sides of the error equations expressed in terms of $\tilde{R}$
lose regularity of only one derivative. Hence, $\tilde{R}$ can be estimated on the desired time span with the help of an appropriate
time--dependent analytic norm which allows one to overcome the loss of regularity. This analytic norm can be used because the approximation function is compactly supported in Fourier space up to an error of order $\mathcal{O}(\epsilon^{3/2})$ in $(H^{s}(\R,\R))^2$.

For infinite depth of water, the NLS approximation can be derived similarly with the difference that the linear dispersion relation \eqref{disprel} has to be replaced by \eqref{disprel2}. The NLS approximation in the case of infinite depth and no surface tension
was justified by Totz and Wu \cite{TW11}. In this situation, the elimination of the 
terms of order $\mathcal{O}(\epsilon)$ in the evolutionary system of the error
is possible without loss of regularity with the help of the above mentioned coordinate transform from \cite{Wu09a}.

The justification of the NLS approximation in case of surface tension is much more complicated and subject of current research. For weak surface tension, there are even situations where the NLS approximation fails; see \cite{SSZ15}.
In \cite{DS06, D16, DH16}, the NLS approximation was justified for various model problems. Each model problem and the water wave problem with surface tension share some of the principal difficulties which have to be addressed in a validity proof of the NLS approximation.

\end{document}